\documentclass[11pt]{article}
\usepackage{amsfonts}
\usepackage{color}
\newtheorem{theo}{Theorem}[section]
\newtheorem{lem}[theo]{Lemma}
\newtheorem{cor}[theo]{Corollary}

\newcommand{\mysection}[1]{\section{#1} \setcounter{equation}{0}}
\newcommand{\proof}{{\sc Proof.} \quad}
\newcommand{\proofc}{{\sc Proof} \ }
\newcommand{\be}{\begin{equation} \label}
\newcommand{\ee}{\end{equation}}
\newcommand{\bea}{\begin{eqnarray}\label}
\newcommand{\eea}{\end{eqnarray}}
\newcommand{\bas}{\begin{eqnarray*}}
\newcommand{\eas}{\end{eqnarray*}}
\newcommand{\bit}{\begin{itemize}}
\newcommand{\eit}{\end{itemize}}
\newcommand{\qed}{\hfill$\Box$ \vskip.2cm}
\newcommand{\nn}{\nonumber}
\newcommand{\R}{\mathbb{R}}

\newcommand{\pO}{\partial\Omega}

\newcommand{\hra}{\hookrightarrow}

\newcommand{\abs}{\\[5pt]}

\newcommand{\io}{\int_\Omega}
\newcommand{\bom}{\overline{\Omega}}

\newcommand{\dw}{D_w}
\newcommand{\dz}{D_z}
\newcommand{\ua}{\underline{a}}

\begin{document}
\title{Dampening effect of logistic source in a two-dimensional haptotaxis system
with nonlinear zero-order interaction}
\author{Zhen Chen\footnote{chenzhendhu@163.com}\\
{\small School of Information Science $\&$ Technology, Donghua University,}\\
{\small Shanghai 200051, P.R.~China} }
\date{}
\maketitle
\begin{abstract}
\noindent
This paper deals with the oncolytic virotherapy model
   \bas
        \left\{ \begin{array}{l}
        u_t = \Delta u - \nabla \cdot (u\nabla v)-uz +\mu u(1-u), \\[1mm]
        v_t = - (u+w)v, \\[1mm]
        w_t = \dw \Delta w - w + uz, \\[1mm]
    z_t = \dz \Delta z - z - uz + \beta w,
         \end{array} \right.
         \qquad (\star)
  \eas
in a bounded domain $\Omega\subset \R^2$ with smooth boundary, where
$\mu, \dw, \dz$ and $\beta$ are prescribed positive parameters.\abs
For any given suitably regular initial data, the global existence of
classical solution to the corresponding homogeneous Neumann
initial-boundary problem for a more general model allowing $\mu=0$
was previously verified in [Y. Tao $\&$ M. Winkler, J. Differential
Equations {\bf 268} (2020), 4973-4997]. This work further shows that
whenever $\mu>0$, the above-mentioned global classical solution to
($\star$) is uniformly bounded; and moreover, if $\beta<1$, then the
solution $(u, v, w, z)$ stabilizes to the constant equilibrium $(1,
0, 0, 0)$ in the topology $L^p(\Omega)\times (L^\infty(\Omega))^3$
with any $p>1$ in a large time limit.\abs
\end{abstract}
\newpage
\mysection{Introduction}
Aim at bypassing the obstacle of drug transfer in traditional
chemotherapy (\cite{swabb}, \cite{jain}), oncolytic virotherapy
becomes an alternative treatment for cancer and it has been
undergoing clinical trials (cf. \cite{coffey} and \cite{russell},
for instance). Very recently, in order to explore efficiency of this
novel therapy, Alzahrani et al. in \cite{eftimie} proposed the
following haptotaxis model
 \be{0}
        \left\{ \begin{array}{lcll}
       u_t &=& \Delta u - \nabla \cdot (u\nabla v) - uz +\mu u(1-u),
        & x\in\Omega, \ t>0, \\[1mm]
    v_t &=& - (u+w)v,
        & x\in\Omega, \ t>0, \\[1mm]
    w_t &=& \dw \Delta w - w + uz,
        & x\in\Omega, \ t>0, \\[1mm]
    z_t &=& \dz \Delta z - z - uz + \beta w,
        & x\in\Omega, \ t>0, \\[1mm]
        & & \hspace*{-15mm}
        (\nabla u -u\nabla v)\cdot \nu=\frac{\partial w}{\partial\nu} = \frac{\partial z}{\partial\nu}=0,
        & x\in\pO, \ t>0, \\[1mm]
    & & \hspace*{-15mm}
        u(x,0)=u_0(x),
        \quad v(x,0)=v_0(x),
        \quad w(x,0)=w_0(x),
    \quad z(x,0)=z_0(x),
        & x\in\Omega,
        \end{array} \right.
\ee with positive parameters $\mu, \beta, \dw$ and $\dz$. Throughout
this work we shall assume that $\Omega\subset\R^2$ is a bounded
domain with smooth boundary. Here $u, w, z$ and $v$ represent the
densities of uninfected tumor cells, infected tumor cells, virus
particles and  normal tissue, respectively. Besides random motion,
the uninfected tumor cells direct their movement toward the higher
densities of tissue, they are assumed to follow a logistic growth
and their number could be reduced due to infection by virus, while
the infected cells randomly diffuse and have a natural death; the
virus is released by infected cells and it experiences diffusion and
loss due to decay and infection; and the tissue is degraded upon
contact with cancer cells. The readers may refer to \cite{eftimie}
for more detailed biological backgrounds behind (\ref{0}).\abs
Different from standard reaction-diffusion equations, the
cross-diffusion term in chemotactic or haptotactic systems like
(\ref{0}) may exhibit a destabilizing feature
(\cite{herrero_velazquez}, \cite{win_jmpa}). Previous qualitative
studies on such types of cross-diffusion systems mainly concentrate
on global existence (\cite{walker_webb}, \cite{taowin_JDE2014},
\cite{zhigun_surulescu_uatay}, \cite{surulescu_win_CMS},
\cite{pang_wang_m3as2018}) and only a few address the large time
behavior of solutions (\cite{litcanu_cmr}, \cite{taowin_sima2015},
\cite{win_jmpa2018}). \abs
In contrast to preceding haptotaxis systems, (\ref{0}) contains a
nonlinear zero-order interaction term $uz$ in the third equation
that gives rise to a new challenge in analysis of this system.
Although (\ref{0}) is indeed globally well-posed
(\cite{taowin_261}), a novel critical parameter phenomenon for
infinite-time blow-up was recently detected for (\ref{0}) without
any growth or sink term (\cite{taowin_262}). Moreover, the global
existence and large time behavior of classical solutions was also
discussed for (\ref{0}) provided that  $\beta<1$ and $\mu=0$
(\cite{taowin_263}). \abs
The purpose of this work is to identify the stabilization effect of
the logistic dampening term in (\ref{0}). To this end, we assume
that
%
%
%
%
%
\be{init}
        \left\{ \begin{array}{l}
        \mbox{$u_0, v_0$ and $w_0$ are nonnegative functions from $C^{2+\vartheta}(\bom)$ for some $\vartheta>0$,} \\[1mm]
        \mbox{with $u_0> 0$, $w_0\not\equiv 0, z_0\not\equiv 0$, $\sqrt{v_0} \in W^{1,2}(\Omega)$ and
    $\frac{\partial u_0}{\partial\nu}=\frac{\partial v_0}{\partial\nu}=\frac{\partial w_0}{\partial\nu}=0$ on $\pO$.}
    \end{array} \right.
\ee We shall prove that whenever $\mu>0$, for any $\beta>0$ the
classical solution of (\ref{0}) is uniformly bounded; furthermore,
it is claimed that if $\beta<1$, then the solution $(u, v, w,z)$
stabilizes to the constant stationary solution $(1, 0, 0, 0)$ in
$L^\infty$ topology. More precisely, we have the following:
\begin{theo}\label{theo1}
  Let $\mu>0$, and let $\Omega\subset\R^2$ be a smoothly bounded domain.
  Then for any given $(u_0,v_0,w_0,z_0)$ satisfying (\ref{init}), the solution $(u,v,w,z)$ of (\ref{0}) is bounded in the sense that
  \be{B}
    \sup_{t>0} \Big\{
    \|u(\cdot,t)\|_{L^\infty(\Omega)} +
    \|v(\cdot,t)\|_{L^\infty(\Omega)} +
    \|w(\cdot,t)\|_{L^\infty(\Omega)} +
    \|z(\cdot,t)\|_{L^\infty(\Omega)} \Big\}
    <\infty.
  \ee
If
  \bas
  \beta \in (0, 1),
  \eas
  then furthermore
  \bas
    u(\cdot,t) \to 1
    \quad \mbox{in } L^p(\Omega)
    \qquad \mbox{for all } p\ge 1
  \eas
  and
  \bas
    v(\cdot,t) \to 0
    \qquad \mbox{in } L^\infty(\Omega)
  \eas
  as well as
  \bas
    w(\cdot,t) \to 0
    \qquad \mbox{in } L^\infty(\Omega)
  \eas
  and
  \bas
    z(\cdot,t) \to 0
    \qquad \mbox{in } L^\infty(\Omega)
  \eas
  as $t\to\infty$.\\
\end{theo}
%
%
%
%
%
%
%
%
%
%
%
%
%
%
%
%
%
%
%
%
%
%
%
%
%
%
%
%
%
%
\mysection{$L^p$ bounds for $w$}
Let us begin with the global smooth solvability of (\ref{0})
previously asserted in \cite{taowin_261}, together with one basic
solution property.
\begin{lem}\label{lem_basic}
  Let $\Omega\subset\R^2$ be a bounded domain with smooth boundary, and assume that
  $(u_0,v_0,w_0,z_0)$ fulfills (\ref{init}). Then the problem (\ref{0}) admits a uniquely determined classical solution
  $(u,v,w,z) \in (C^{2,1}(\bom\times [0,\infty)))^4$ for which $v$ is nonnegative, and for which $u,w$ and $z$ are positive
  in $\bom\times (0,\infty)$. Moreover,
   \be{vinfty}
    \|v(\cdot,t)\|_{L^\infty(\Omega)} \le \|v_0(\cdot)\|_{L^\infty(\Omega)}
    \qquad \mbox{for all } t>0.
  \ee
%
\end{lem}
The presence of logistic dampening in the first equation in
(\ref{0}) will be decisive for us to further verify the uniform
boundedness of the global solution constructed in Lemma
\ref{lem_basic}, and correspondingly the assumption that $\mu>0$
will play a key role in deriving the following $L^1$ bounds for $u,
w$ and $z$ and a space-time $L^2$ bound for $u$, which serve as a
starting point of our subsequently reasoning.
\begin{lem}\label{lem3}
 There exists $C>0$ such that
  \be{3.1}
    \io u(\cdot, t)\le C,
    \quad \io w(\cdot, t)\le C
    \quad \mbox{and}\quad
    \io z(\cdot, t)\le C
      \ee
  as well as
  \be{3.2}
    \int_t^{t+1} \io u^2 \le C
  \ee
 for all $t>0$.
\end{lem}
\proof
  We firstly assume that $\beta>0$ and use (\ref{0}) to compute
  \bas
    \frac{d}{dt} \bigg\{ 2\beta\io u +2\beta\io w + \io z \bigg\}
    &=& \bigg\{ -2\beta \io uz + 2\beta\mu \io u -2\beta \mu \io u^2 \bigg\} \\
    & &  + \bigg\{ - 2\beta \io w + 2\beta \io uz \bigg\} \\
    & &  + \bigg\{ -\io z -\io uz + \beta \io w \bigg\}\\
    &=& 2\beta\mu \io u -\beta \io w -\io z -\io uz -2\beta \mu \io u^2
  \eas
for all $t>0$, so that by Young's inequality,
  \bas
    \frac{d}{dt}  \bigg\{ 2\beta\io u +2\beta\io w + \io z \bigg\}
   &+& \frac{1}{2} \cdot \bigg\{ 2\beta\io u +2\beta\io w + \io z \bigg\}
    + \beta \mu \io u^2\\
   &= & - \beta (1+2\mu)\io u -\beta \mu \io u^2 -\frac{1}{2}\io z -\io uz\\
   &\le& - \beta (1+2\mu)\io u -\beta \mu \io u^2\\
   &\le& \frac{\beta (1+2\mu)^2}{4\mu}\cdot |\Omega|
    \qquad \mbox{for all } t>0.
  \eas
  Upon an ODE comparison, this readily implies (\ref{3.1}) due to $\beta>0$ and thereby a further integration
  over time yields (\ref{3.2}) thanks to the fact that $\mu>0$. If
  $\beta=0$, by estimating $\frac{d}{dt}\Big\{\io w +\io z\Big\}$ and
  $\frac{d}{dt}\io u$ we still can arrive at (\ref{3.1}) and
  (\ref{3.2}) in a same manner.
 \qed
With the help of parabolic smooth properties in the spatial
two-dimensional setting, the $L^1$ bound of $w$ actually implies
$L^p$ integrability of $z$ for arbitrarily large finite p.
\begin{lem}\label{lem4}
  Let $p\in (1,\infty)$. Then there exists $C(p)>0$ such that
  \be{4.1}
    \|z(\cdot,t)\|_{L^p(\Omega)} \le C(p)
    \qquad \mbox{for all } t>0.
  \ee
\end{lem}
\proof
   By nonnegativity of $uz$ and $z$, in light of the order preserving property of $e^{\sigma\Delta}$ for $\sigma\ge 0$ we firstly
  find that
  \bas
    z(\cdot,t)
    &=& e^{t(\dz\Delta-1)} z_0
    - \int_0^t e^{(t-s)(\dz\Delta-1)} u(\cdot,s)z(\cdot,s) ds
    + \beta \int_0^t e^{(t-s)(\dz\Delta-1)} w(\cdot,s) ds \\
    &\le& e^{t\dz\Delta} z_0
    + \beta \int_0^t e^{(t-s)(\dz\Delta-1)} w(\cdot,s) ds
    \quad \mbox{in } \Omega
    \qquad \mbox{for all } t>0,
  \eas
  and hence, by well-known smoothing properties of the Neumann heat semigroup $(e^{\sigma\Delta})_{\sigma\ge 0}$
  (\cite{win_JDE2010}) and Lemma \ref{lem3}, secondly see that
  \bas
    \|z(\cdot,t)\|_{L^p(\Omega)}
    &\le& \|e^{t\dz\Delta} z_0\|_{L^p(\Omega)}
    + c_1\beta \int_0^t \Big(1+(t-s)^{-1+\frac{1}{p}}\Big) e^{-(t-s)} \|w(\cdot,s)\|_{L^1(\Omega)} ds \\
    &\le& \|z_0\|_{L^p(\Omega)}
    + c_1 c_2 \beta \int_0^t \Big(1+(t-s)^{-1+\frac{1}{p}}\Big) e^{-(t-s)} ds
    \qquad \mbox{for all } t>0
  \eas
  with some $c_1>0$ and $c_2>0$. This entails (\ref{4.1}), because
   \bas
    \int_0^t \Big(1+(t-s)^{-1+\frac{1}{p}}\Big) e^{-(t-s)} ds
    \le \int_0^\infty (1+\sigma^{-1+\frac{1}{p}}) e^{-\sigma}
    d\sigma < 2+p
    \qquad \mbox{for all } t>0.
  \eas
  \qed
In light of the above $L^p$ integrability information of $z$ along
with a space-time $L^2$ bound for $u$ in (\ref{3.2}) and according
to parabolic smooth estimates in the two-dimensional case once
again, we also obtain $L^p$ bounds of $w$ for arbitrary $p>2$.
\begin{lem}\label{lem5}
Let $p\in (2, \infty)$. Then there exists $C(p)>0$ with the property
  \be{5.1}
    \|w(\cdot,t)\|_{L^p(\Omega)} \le C(p)
    \qquad \mbox{for all } t>0.
  \ee
\end{lem}
\proof Given $p>2$, we choose $q:=\frac{2p}{p+1}\in (1, 2)$.
   We then recall known smoothing estimates for the Neumann heat semigroup on $\Omega$ (\cite{win_JDE2010}) to obtain positive constants
  $c_1=c_1(p)$ and $c_2=c_2(p)$ satisfying
  \be{5.2}
    \|e^{\dw \Delta} \varphi\|_{L^p(\Omega)} \le c_1 \|\varphi\|_{L^1(\Omega)}
    \qquad \mbox{for all } \varphi\in C^0(\bom)
  \ee
  and
  \be{5.3}
    \|e^{\sigma\dw \Delta}\varphi\|_{L^p(\Omega)} \le c_2 \sigma^{-\frac{p-1}{2p}} \|\varphi\|_{L^q(\Omega)}
    \qquad \mbox{for all $\sigma\in (0,1)$ and each } \varphi\in
    C^0(\bom).
  \ee
  Now in light of a variation-of-constants representation of $w$ solving the third equation in (\ref{0}), we can invoke (\ref{5.2})
  and (\ref{5.3}) to estimate
  \bea{5.4}
    \|w(\cdot,t)\|_{L^p(\Omega)}
    &=& \bigg\| e^{\dw\Delta -1} w(\cdot,t-1) + \int_{t-1}^t e^{(t-s)(\dw\Delta-1)} u(\cdot,s) z(\cdot,s) ds \bigg\|_{L^p(\Omega)}
        \nn\\
    &\le& e^{-1} \cdot c_1 \|w(\cdot,t-1)\|_{L^1(\Omega)}
    + c_2\int_{t-1}^t (t-s)^{-\frac{p-1}{2p}} \|u(\cdot,s)z(\cdot,s)\|_{L^q(\Omega)} ds
  \eea
  for all $t>1$,
  where according to the H\"older inequality and Lemma \ref{lem4}, we find $c_3=c_3(p)>0$
  such that
  \bas
    c_2\int_{t-1}^t (t-s)^{-\frac{p-1}{2p}} \|u(\cdot,s)z(\cdot,s)\|_{L^q(\Omega)} ds
    &\le& c_2 \int_{t-1}^t (t-s)^{-\frac{p-1}{2p}} \|u(\cdot,s)\|_{L^2(\Omega)} \|z(\cdot,s)\|_{L^\frac{2q}{2-q}(\Omega)} ds \\
    &\le& c_3\int_{t-1}^t (t-s)^{-\frac{p-1}{2p}} \|u(\cdot,s)\|_{L^2(\Omega)} ds \\
    &\le& c_3 \cdot \bigg\{ \int_{t-1}^t (t-s)^{-\frac{p-1}{p}} ds \bigg\}^\frac{1}{2} \cdot
        \bigg\{ \int_{t-1}^t \|u(\cdot,s)\|_{L^2(\Omega)}^2 ds \bigg\}^\frac{1}{2} \\
    &=& p c_3 \cdot
        \bigg\{ \int_{t-1}^t \|u(\cdot,s)\|_{L^2(\Omega)}^2 ds \bigg\}^\frac{1}{2}
    \qquad \mbox{for all } r>1.
  \eas
  Since
  \bas
     \sup_{t>1} \|w(\cdot,t-1)\|_{L^1(\Omega)} <\infty
    \qquad \mbox{and} \qquad
    \int_{t-1}^t \|u(\cdot,s)\|_{L^2(\Omega)}^2 ds
    =\int_{t-1}^t \io u^2 <\infty
    \qquad \mbox{for all } t>1
  \eas by Lemma \ref{lem3}, (\ref{5.4}) yields (\ref{5.1}), because
  $w$ is bounded in $\bom\times [0,1]$ by Lemma \ref{lem_basic}.
\qed
\mysection{$L^\infty$ estimates on $u$}
In this section we shall establish the boundedness of $u$ in
$L^\infty$. For this purpose, as performed in
\cite{fontelos_friedman_hu}, \cite{friedman_tello},
\cite{walker_webb} and \cite{taowin_JDE2014} which dealt with
haptotaxis-related systems,  we introduce the variable change
 \be{a}
    a:= u e^{-v}
\ee and then in view of (\ref{0}),
 \be{0a}
        \left\{ \begin{array}{lcll}
        a_t &=& e^{-v} \nabla \cdot (e^v \nabla a)
    + a(ae^v+w)v - az +\mu a (1-ae^v),
        & x\in\Omega, \ t>0, \\[1mm]
        v_t &=& -(ae^v+w)v,
        & x\in\Omega, \ t>0, \\[1mm]
        & & \hspace*{-13mm}
        \frac{\partial a}{\partial\nu}=0,
        & x\in\pO, \ t>0, \\[1mm]
        & & \hspace*{-13mm}
        a(x,0)=u_0(x) e^{-v_0(x)} =: a_0(x),
        \quad v(x,0)=v_0(x),
          & x\in\Omega.
        \end{array} \right.
\ee
A direct testing procedure leads to the following primary inequality
which will be invoked in Lemma \ref{lem6} and Lemma \ref{lem8}
below.
\begin{lem}\label{lem2}
  Let $p>1$. Then
  \be{1.1}
    \frac{d}{dt} \io e^v a^p
    \le - p(p-1)\io e^v a^{p-2} |\nabla a|^2
    + (p-1) \io e^v a^p (ae^v +w)v
    +\mu p \io e^v a^p
    \qquad \mbox{for all } t>0.
  \ee
\end{lem}
\proof
  From (\ref{0a}) we infer that $(ae^v+w)v=-v_t$ in $\Omega\times (0,\infty)$
  and that
  \bas
    \frac{d}{dt} \io e^v a^p
    &=& p \io e^v a^{p-1} \cdot \Big\{ e^{-v} \nabla\cdot (e^v \nabla a) - av_t - az +\mu a -\mu a^2 e^v\Big\}
    + \io e^v a^p v_t\\
    &=& -p(p-1)\io e^v a^{p-2} |\nabla a|^2
    - (p-1) \io e^v a^p v_t  +\mu p \io e^v a^p \\
    & & -\rho p \io e^v a^p z
        -\mu p\io e^{2v} a^{p+1}
  \eas
for all $t>0$, this leads to (\ref{1.1}) upon abandoning the
rightmost two nonpositive summands. \qed
Using a bootstrap $L^p$-estimate technique developed in \cite[Lemma
3.8]{taowin_zamp2016} and relying on (\ref{3.2}), Lemma \ref{lem5}
and Lemma \ref{lem_basic}, an application of Lemma \ref{lem2} to
$p=2$ can yield a bound for $a$ in $L^2$.
\begin{lem}\label{lem6}
  There exists $C>0$ fulfilling
  \be{6.1}
    \|a(\cdot,t)\|_{L^2(\Omega)} \le C
    \qquad \mbox{for all } t>0.
  \ee
\end{lem}
\proof Lemma \ref{lem_basic} guarantees the existence of $c_1>0$
such that
 \be{6.2}
 \|v(\cdot, t)\|_{L^\infty(\Omega)} \le c_1
 \qquad\mbox{for all } t>0,
 \ee
whereas Lemma \ref{lem5} provides $c_2>0$ satisfying
  \be{6.3}
    \io w^3(\cdot,t) \le c_2
    \qquad \mbox{for all } t>0.
  \ee
  Now from Lemma \ref{lem2} we infer the inequality
  \be{6.4}
    \frac{d}{dt} \io e^v a^2
    + 2\io e^v |\nabla a|^2
    \le \io e^{2v} a^3 v
    + \io e^v a^2 vw
    + 2\mu \io e^v a^2
    \qquad \mbox{for all } t>0,
  \ee
  where by Young's inequality, (\ref{6.2}) and (\ref{6.3}),
  \bas
  \hspace*{-4mm}
    \io e^{2v} a^3 v
    + \io e^v a^2 vw
    + 2\mu \io e^v a^2
    &\le&  c_1e^{2c_1} \io  a^3
    + c_1 e^{c_1}\io a^2 w
    + 2\mu e^{c_1} \io a^2 \\
    &\le& \Big\{c_1e^{2c_1} +c_1e^{c_1} + 2\mu e^{c_1}\Big\} \cdot \io a^3 + c_1 e^{c_1} \io w^3 + 2\mu e^{c_1} \cdot |\Omega|\\
    &\le& c_3 \io a^3 + c_4
    \qquad \mbox{for all } t>0
  \eas
  with $c_3:= c_1e^{2c_1} +c_1e^{c_1} + 2\mu e^{c_1}$ and $c_4:= c_1 e^{c_1} c_2 + 2\mu e^{c_1} \cdot |\Omega|$.
  As the Gagliardo-Nirenberg inequality along with Young's inequality yields that
  \bas
  c_3 \io a^3
  = c_3\|a\|_{L^3(\Omega)}^3
  &\le& c_5\|\nabla a\|_{L^2(\Omega)} \|a\|_{L^2(\Omega)}^2 +
         c_5\|a\|_{L^2(\Omega)}^3\\
  &\le&  \Big\{2\|\nabla a\|_{L^2(\Omega)}^2
       +c_5^2\|a\|_{L^2(\Omega)}^4\Big\}
      +\Big\{\|a\|_{L^2(\Omega)}^4 +c_5^4\Big\}\\
  &\le& 2\io e^v |\nabla a|^2 +c_6\|a\|_{L^2(\Omega)}^4 +c_6
    \qquad \mbox{for all } t>0
  \eas
 with some $c_5>0$ and $c_6:=\max\{c_5^2+1,\, c_5^4\}$, by$\io a^2\le \io e^v a^2$
 we thus infer from (\ref{6.4}) that
  \be{6.5}
    \frac{d}{dt} \io e^v a^2
    \le \bigg\{c_6\io a^2(\cdot, t)\bigg\}\cdot \io e^v a^2 +c_6
    \qquad \mbox{for all } t>0.
  \ee Letting
  \bas
 h(t):=c_6\io a^2(\cdot, t)
 \qquad\mbox{for all } t>0,
  \eas
 We may apply Lemma \ref{lem3} along with (\ref{vinfty}) to find
 $c_7>0$ fulfilling
  \be{6.6}
  \int_t^{t+1} h(s) ds \le c_7
  \qquad\mbox{for all } t>0,
  \ee
  and thereby for any given $t>1$, it is possible to fix
  $t_0=t_0(t)\in (t-1, t)$ such that $t_0>0$ and
  \bas
 \io a^2(\cdot, t_0)\le c_8:= \max\bigg\{\io a_0^2, \, c_7\bigg\},
  \eas
  so that  integrating (\ref{6.5}) shows that
  \bas
 \io e^{v(\cdot, t)} a^2(\cdot, t)
 &\le& \bigg(\io e^{v(\cdot, t_0)} a^2(\cdot, t_0) \bigg) \cdot
     e^{\int_{t_0}^t h(s) ds}
     + c_6\int_{t_0}^t e^{\int_s^t h(\sigma) d\sigma} ds\\
 &\le& e^{c_1}c_8\cdot e^{c_7} + c_6 e^{c_7}
 \qquad\mbox{for all } t>1
  \eas
  due to (\ref{6.2}), (\ref{6.6}) and the fact that $t-t_0\le 1$.
  Since $e^v\ge 1$ and $a$ is bounded in $\bom\times [0,1]$ by Lemma
  \ref{lem_basic},
  this proves (\ref{6.1}).
  \qed
Using this $L^2$ bound of $a$ together with the smoothing properties
of the Neumann heat semigroup on $\Omega$ once more, we can further
improve the regularity of $w$ and $z$ in the following sense.
\begin{lem}\label{lem7}
  There exists $C>0$  with the property that
  \be{7.1}
    \|w(\cdot,t)\|_{W^{1, 3}(\Omega)} \le C
    \qquad \mbox{for all } t>0
  \ee
  and
  \be{7.2}
    \|z(\cdot,t)\|_{W^{1, 3}(\Omega)} \le C
    \qquad \mbox{for all } t>0.
  \ee
\end{lem}
\proof Combine Lemma \ref{lem6}, Lemma \ref{lem5} and Lemma
\ref{lem4} with (\ref{a}) and (\ref{vinfty}) we can fix positive
constants $c_1, c_2$ and $c_3$ satisfying
  \bas
    \|u(\cdot,t)\|_{L^2(\Omega)} \le c_1,
    \quad
    \|w(\cdot,t)\|_{L^2(\Omega)} \le c_2
    \quad \mbox{and} \quad
    \|z(\cdot,t)\|_{L^4(\Omega)} \le c_3
    \qquad \mbox{for all } t>0.
  \eas
  Since well-known smoothing properties of the Neumann heat semigroup on $\Omega$ (\cite{win_JDE2010}) ensure
  the existence of $c_4>0$ and $c_5>0$ such that
  \bas
    \|e^{\dw\Delta} \varphi\|_{W^{1,3}(\Omega)} \le c_4\|\varphi\|_{L^2(\Omega)}
    \qquad \mbox{for all } \varphi\in C^0(\bom)
  \eas
  and
  \bas
    \|e^{\sigma\dw\Delta}\varphi\|_{W^{1,3}(\Omega)} \le c_5 \sigma^{-\frac{11}{12}} \|\varphi\|_{L^\frac{4}{3}(\Omega)}
    \qquad \mbox{for all $\sigma\in (0,1)$ and any } \varphi\in
    C^0(\bom),
  \eas we thus see that by the H\"{o}lder inequality,
  \bas
    \|w(\cdot,t)\|_{W^{1,3}(\Omega)}
    &=& \bigg\| e^{\dw\Delta-1} w(\cdot,t-1) + \int_{t-1}^t e^{(t-s)(\dw\Delta-1)} u(\cdot,s)z(\cdot,s) ds \bigg\|_{W^{1,3}(\Omega)}
        \\
    &\le& c_4 e^{-1} \|w(\cdot,t-1)\|_{L^2(\Omega)}
    + c_5 \int_{t-1}^t (t-s)^{-\frac{11}{12}} \|u(\cdot,s)z(\cdot,s)\|_{L^\frac{4}{3}(\Omega)} ds \\
    &\le& c_4 e^{-1} \|w(\cdot,t-1)\|_{L^2(\Omega)}
    + c_5 \int_{t-1}^t (t-s)^{-\frac{11}{12}} \|u(\cdot,s)\|_{L^2(\Omega)} \|z(\cdot,s)\|_{L^4(\Omega)} ds \\
    &\le& c_2 c_4 e^{-1}
    + 12c_1 c_3 c_5
    \qquad \mbox{for all } t>1,
  \eas which entails (\ref{7.1}), while (\ref{7.2}) can be proved in
  a quite similar manner. \qed
Lemma \ref{lem7} actually implies the boundedness of $w$ in
$L^\infty$ due to the Sobolev embedding $W^{1,p}(\Omega) \hra
C^0(\bom)$ for $p>n$, with $n$ denoting the space dimension.
Depending on this and Lemma \ref{lem6}, we can achieve the
boundedness of $a$ in $L^\infty$ via a Moser iteration technique.
\begin{lem}\label{lem8}
  There exists $C>0$ such that
  \be{8.1}
    \|a(\cdot,t)\|_{L^\infty(\Omega)} \le C
    \qquad \mbox{for all } t>0.
  \ee
\end{lem}
\proof According to Lemma \ref{lem_basic} and Lemma \ref{lem7} along
with the Sobolev embedding $W^{1,3}(\Omega) \hra C^0(\bom)$, we can
pick $c_1>0$ and $c_2>0$ such that $v\le c_1$ and $w\le c_2$ on
$\Omega\times (0, \infty)$. From Lemma \ref{lem2} along with Young's
inequality we next can infer that for all $p>1$,
  \be{8.2}
    \frac{d}{dt} \io e^v a^p
    + 3\io |\nabla a^\frac{p}{2}|^2
    + \io e^v a^p
    \le c_3 p \io a^{p+1}
    + c_4
    \qquad \mbox{for all } t>0
  \ee with $c_3>0$ and $c_4>0$ which are independent of $p$. Since
Lemma \ref{lem6} provides $c_5>0$ such that
 \bas
 \io a^2(\cdot, t) \le c_5
 \qquad\mbox{for all } t>0,
 \eas
we thus can invoke a Moser-type iteration method to derive
(\ref{8.1}) (cf. \cite[Lemma 6.4]{taowin_263} for details). \qed
\mysection{Stabilization when $\beta<1$}
Throughout this section we shall suppose that $\beta<1$, and we
shall investigate the asymptotic behavior of the global classical
solution constructed in Lemma \ref{lem_basic}. In fact, the
assumption that $\beta<1$ readily results in the exponential decay
property of $w$ and $z$.
\begin{lem}\label{lem9}
 Let $\beta<1$. Then there exist $\gamma \in (0,1)$ and $C>0$ such
 that
  \be{9.1}
    \|w(\cdot,t)\|_{L^\infty(\Omega)} \le C e^{-\gamma t}
    \qquad \mbox{for all } t>0
  \ee
  and
  \be{9.2}
    \|z(\cdot,t)\|_{L^\infty(\Omega)} \le C e^{-\gamma t}
    \qquad \mbox{for all } t>0
  \ee
\end{lem}
\proof We use the third and fourth equations in (\ref{0}) to compute
  \bas
    \frac{d}{dt} \bigg\{ \io w + \io z \bigg\}
    &=& \bigg\{ - \io w + \io uz \bigg\}
    + \bigg\{ - \io z - \io uz + \beta \io w \bigg\} \\
    &=& - (1-\beta) \io w - \io z \\
    &\le& -(1-\beta) \cdot \bigg\{ \io w + \io z \bigg\}
    \qquad \mbox{for all } t>0
  \eas
  which implies that
  \bas
    \io w + \io z \le \bigg\{ \io w_0 + \io z_0 \bigg\} \cdot e^{-(1-\beta)t}
    \qquad \mbox{for all } t>0.
  \eas
This in conjunction with Lemma the Gagliardo-Nirenberg inequality
and Lemma \ref{lem7} yields $c_1>0$ and $c_2>0$ such that
 \bas
    \|w(\cdot,t)\|_{L^\infty(\Omega)}
    &\le& c_1\|w(\cdot,t)\|_{W^{1,3}(\Omega)}^\frac{6}{7}
         \|w(\cdot,t)\|_{L^1(\Omega)}^\frac{1}{7}\\
    &\le& c_2 e^{-\frac{1-\beta}{7}t}
    \qquad \mbox{for all } t>0,
  \eas
and thus the proof is complete with $\gamma:=\frac{1-\beta}{7}$ and
$C:=c_2$. \qed
Since $z$ eventually vanishes in $L^\infty$ due to (\ref{9.2}), from
(\ref{0a}) and the parabolic comparison principle we can infer that
$a$ has a positive lower bound, and thus so does $u$ in view of
(\ref{a}) and (\ref{vinfty}).
\begin{lem}\label{lem10}
If $\beta<1$, then  one can find $\delta>0$  fulfilling
  \be{10.1}
   u(x,t)\ge \delta
    \qquad \mbox{for all } (x, t)\in \Omega\times (0, \infty).
  \ee
\end{lem}
\proof According to Lemma \ref{lem9} and Lemma \ref{lem_basic}, it
is possible for us to pick $t_0>0$ and $c_1>0$ such that
 \bas
  z(x, t)\le c_1 e^{-\gamma t} \le \frac{\mu}{2}
   \qquad \mbox{for all $x\in\Omega$ and } t>t_0.
 \eas
From the first equation in (\ref{0a}) along with this and Lemma
\ref{lem_basic} once again we infer that $a$ satisfies
 \bas
    a_t \ge e^{-v} \nabla \cdot (e^v \nabla a)  +\mu a
    \Big(\frac{1}{2} -e^{\|v_0\|_{L^\infty(\Omega)}} a\Big)
    \qquad \mbox{for all $x\in\Omega$ and } t>t_0
  \eas
  with $\frac{\partial a}{\partial\nu}=0$ throughout $\pO\times  (t_0,\infty)$.
  On the other hand, for
  \bas
    \ua(x,t):=y(t),
    \qquad x\in\bom, \ t\ge t_0,
  \eas
  with $y\in C^1([t_0,\infty))$ denoting the solution of
  \bas
    \left\{ \begin{array}{l}
   y'(t) = \mu y \Big(\frac{1}{2}-c_2 y\Big),
    \qquad t>t_0, \\[1mm]
    y(t_0):=\inf_{x\in\Omega}a(x, t_0)>0,
    \end{array} \right.
  \eas where $c_2:=e^{\|v_0\|_{L^\infty(\Omega)}}$. Then
  we have
  \bas
    & & \hspace*{-20mm}
    \ua_t - e^{-v(x,t)} \nabla \cdot (e^{v(x,t)} \nabla \ua) -\mu \ua \Big(\frac{1}{2}-c_2\ua\Big)\\
    &=& y'(t) - \mu y \Big(\frac{1}{2}-c_2 y\Big) \\[1mm]
    &=& 0
    \qquad \mbox{for all $x\in\Omega$ and } t>t_0,
  \eas
  and since clearly $\frac{\partial \ua}{\partial\nu}=0$ on $\pO\times (t_0,\infty)$ and
  \bas
    \ua(x,t_0)= y(t_0) = \inf_{x\in\Omega} a(x,t_0) \le a(x,t_0)
    \qquad \mbox{for all } x\in\Omega,
  \eas
 an application of the comparison principle readily shows that $a\ge \ua$ in $\Omega\times (t_0,\infty)$ and hence
  \be{10.2}
    u(x,t) = a(x,t) e^{v(x,t)}
    \ge a(x,t) \ge y(t)
    \qquad \mbox{for all $x\in\Omega$ and } t>t_0.
  \ee
Since we have $y(t)\to \frac{1}{2c_2}$ as $t\to\infty$ due to
$y(t_0)>0$, we can conclude that there exists $t_1\ge t_0$ such that
 \be{10.3}
 y(t)\to \frac{1}{4c_2}
 \qquad\mbox{for all } t>t_1.
 \ee
 Finally, the positivity and continuity of $u(x, t)$ in $\Omega
 \times (0, \infty)$ by Lemma \ref{lem_basic} guarantees the
 existence of $c_3>0$ fulfilling
 \be{10.4}
 c_3:=\min_{\bom \times [0, t_1]} u(x, t)>0.
 \ee
So that (\ref{10.1}) results from (\ref{10.2})-(\ref{10.4}) with
$\delta :=\min\Big\{\frac{1}{4c_2}, \, c_3\Big\}>0$.
 \qed
According to the second equation in (\ref{0}), Lemma \ref{lem10}
along with Lemma \ref{lem_basic} immediately leads to the following.
\begin{cor}\label{cor11}
Let $\beta<1$. Then we have
  \be{11.1}
    v(x,t) \le \|v_0(\cdot)\|_{L^\infty(\Omega)} \cdot e^{-\delta t}
    \qquad  \mbox{for all $x\in\Omega$ and } t>0,
   \ee
   where $\delta>0$ is defined by Lemma \ref{lem10}.
\end{cor}
\proof
  From (\ref{0}), Lemma \ref{lem10} and the nonnegativity of $v$ and $w$ we obtain that
  \bas
    v_t = -(u+w)v \le -\delta v
    \qquad\mbox{for all $x\in\Omega$ and } t>0,
  \eas
 from which (\ref{11.1}) readily follows.
\qed
With the outcomes of Lemma \ref{lem8}, Lemma \ref{lem9} and
Corollary \ref{cor11}, we firstly obtain the following weak
convergence information of $a$, which is based on a specific testing
procedure.
\begin{lem}\label{lem12}
If $\beta<1$, then  there holds
  \be{12.1}
   \int_0^\infty \io (ae^v-1)^2<\infty.
  \ee
\end{lem}
\proof Clearly, Lemma \ref{lem8}, Lemma \ref{lem9} and Corollary
\ref{cor11} ensure the existence of $c_1>0$, $c_2>0$, $c_3>0$,
$c_4>0$ and $c_5>0$ such that
 \be{12.2}
a\le c_1, \quad v\le c_2, \quad\mbox{and}\quad
 w\le c_3
 \qquad\mbox{in } \Omega\times (0, \infty)
 \ee
 as well as
 \be{12.3}
\int_0^\infty \io v \le c_4 \quad\mbox{and}\quad
 \int_0^\infty \io z \le c_5.
 \ee
According to the basic inequality: $s-1-\ln s\ge 0$ for all $s> 0$
and in view of the fact that $v_t\le 0$, we use (\ref{0a}) to
estimate
 \bas
 \frac{d}{dt}\io e^v(a -1 -\ln a)
 &=& \io e^v \Big(1-\frac{1}{a} \Big) a_t +\io e^v (a-1-\ln a ) v_t\\
 &\le& \io e^v \Big(1-\frac{1}{a} \Big) a_t\\
 &=& \io e^v \Big(1-\frac{1}{a} \Big) \cdot \Big\{ e^{-v}\nabla
 \cdot (e^v \nabla a) +a (a e^v +w)v -az +\mu a (1-ae^v)\Big\}\\
 &=& -\io e^v \frac{|\nabla a|^2}{a^2}
     +\io e^v (a-1) (a e^v +w)v -\io e^v (a-1) z \\
 & & +\mu \io (ae^v-e^v)(1-ae^v)\\
 &\le& \io e^v a (a e^v +w)v +\io e^v z +\mu \io(1-e^v)(1-ae^v)\\
 &  &  -\mu \io (ae^v -1)^2
 \qquad\mbox{for all } t>0
 \eas
due to the nonnegativity of $a, v, w$ and $z$. Here by the
elementary inequality: $1-e^{-s}\le s$ for all $s\ge 0$, we have
$(e^v-1)=e^v(1-e^{-v})\le e^v v$, and thus in light of this and
(\ref{12.2}) we see that
 \bas
 \hspace*{-6mm}
 \io e^v a (a e^v +w)v +\io e^v z &+& \mu \io(1-e^v)(1-ae^v)\\
 &\le& \io e^v a (a e^v +w)v +\io e^v z +\mu \io e^v(ae^v- 1)v\\
 &\le& \io \Big(a^2e^{2v} +aw e^v +\mu a e^{2v}\Big) v +\io e^v z\\
 &\le& \Big( c_1^2 e^{2c_2} +c_1 c_3 e^{c_2} +\mu c_1 e^{2c_2}\Big)
      \io v +e^{c_2} \io z
      \qquad \mbox{for all } t>0.
 \eas
Therefore, we arrive at
 \bas
 \frac{d}{dt}\io e^v(a -1 -\ln a) + \mu \io (ae^v -1)^2 \le c_6 \io
 v+c_7 \io z
 \qquad \mbox{for all } t>0
 \eas
with $c_6:=c_1^2 e^{2c_2} +c_1 c_3 e^{c_2} +\mu c_1 e^{2c_2}$ and
$c_7:=e^{c_2}$. Integrating this in time entails that
 \bas
 \hspace*{-5mm}
\io e^{v(\cdot, t)}\{a(\cdot, t) -1 -\ln a(\cdot, t)\}
  &+&    \mu \int_0^t \io (ae^v -1)^2 \\
  &\le&   \io e^{v_0} (a_0 -1 -\ln a_0) + c_6 \int_0^t\io v+c_7 \int_0^t\io z\\
  &\le&   \io e^{v_0} (a_0 -1 -\ln a_0) + c_6 \int_0^\infty\io v+c_7 \int_0^\infty\io  z\\
  &=&   \io e^{v_0} (a_0 -1 -\ln a_0) + c_4c_6 +c_5 c_7
  \qquad\mbox{for all } t>0.
 \eas
Since  $s-1-\ln s\ge 0$ for all $s> 0$ as mentioned before and since
$a(x, t)>0$ in $\Omega\times (0, \infty)$ by Lemma \ref{lem_basic},
this implies that
 \bas
 \int_0^\infty \io (ae^v -1)^2 \le c_8 :=\frac{1}{\mu} \cdot \bigg\{\io e^{v_0} (a_0 -1 -\ln a_0) + c_4c_6 +c_5 c_7
     \bigg\}<\infty
 \eas
thanks to $\mu>0$ and the fact that $a_0:=u_0 e^{-v_0}>0$ by
(\ref{init}), and thereby completes the proof.
 \qed
In order to improve the above weak stabilization information, we
need further regularity properties of $a$, which can be obtained
through another testing procedure.
\begin{lem}\label{lem13}
Let $\beta<1$, then
  \be{13.1}
   \int_0^\infty \io a_t^2<\infty
  \ee
  and
  \be{13.2}
 \sup_{t>0} \io |\nabla a(\cdot, t)|^2<\infty.
  \ee
\end{lem}
\proof  Lemma \ref{lem8}, Lemma \ref{lem9} and Corollary \ref{cor11}
warrant the boundedness of $a, v, w$ and $z$ as well as the
exponential decay of $v$ and $z$, which enables us to find $c_i>0,
i\in \{1, \cdots, 5\}$, fulfilling
 \be{13.3}
 2e^v a^2(ae^v +w^2) \le c_1, \quad
 2e^v a^2\le c_2 \quad\mbox{and}\quad
 e^va^2\le c_3
 \qquad\mbox{in } \Omega\times (0, \infty)
 \ee
and moreover
 \be{13.4}
 \int_0^\infty\io v^2\le c_4 \quad\mbox{and}\quad
 \int_0^\infty\io z^2\le c_5.
 \ee
Lemma \ref{lem12} also provides $c_6>0$ satisfying
 \be{13.5}
  \int_0^\infty \io (ae^v-1)^2<c_5.
 \ee
Next, multiplying the first equation in (\ref{0a}) by $e^v a_t$ and
integrating by parts, we have
  \bea{13.6}
    \io e^v a_t^2
    &=& \io a_t \nabla \cdot (e^v\nabla a)  + \io e^v a_t a(ae^v +w)v
        - \io e^v a_t az + \mu \io e^v a_t a(1-ae^v)\nn\\
    &=& - \io e^v \nabla a\cdot\nabla a_t
        + \io e^v a_t a(ae^v +w)v
        - \io e^v a_t az + \mu \io e^v a_t a(1-ae^v)
  \eea
for all $t>0$,  where by nonpositivity of $v_t$,
  \bea{13.7}
    - \io e^v \nabla a\cdot\nabla a_t
    &=& - \frac{1}{2} \io e^v \partial_t |\nabla a|^2 \nn\\
    &=& - \frac{1}{2} \frac{d}{dt} \io e^v |\nabla a|^2
    + \frac{1}{2} \io e^v |\nabla a|^2 v_t \nn\\
    &\le& - \frac{1}{2} \frac{d}{dt} \io e^v |\nabla a|^2
    \qquad \mbox{for all } t>0.
  \eea
Moreover, in view of Young's inequality and (\ref{13.3}),
  \bea{13.8}
 \io e^v a_t a(ae^v +w)v
 &\le& \frac{1}{8} \io e^v a_t^2 + 2\io e^v a^2 (ae^v+w)^2 v^2 \nn\\
 &\le& \frac{1}{8} \io e^v a_t^2 + c_1\io v^2
  \qquad\mbox{for all } t>0
  \eea
and similarly,
  \bea{13.9}
 - \io e^v a_t az
 &\le& \frac{1}{8} \io e^v a_t^2 + 2\io e^v a^2 z^2 \nn\\
 &\le& \frac{1}{8} \io e^v a_t^2 + c_2\io z^2
  \qquad\mbox{for all } t>0
  \eea
as well as
  \bea{13.10}
 \mu \io e^v a_t a(1-ae^v)
 &\le& \frac{1}{4} \io e^v a_t^2 + \mu^2\io e^v a^2 (1-ae^v)^2 \nn\\
 &\le& \frac{1}{4} \io e^v a_t^2 + \mu^2 c_3\io (1-ae^v)^2
  \qquad\mbox{for all } t>0.
  \eea
Collecting (\ref{13.6})-(\ref{13.10}) and integrating in time we
obtain
 \bas
 \hspace*{-5mm}
  \frac{1}{2}\int_0^t\io e^v a_t^2
  &+ &\frac{1}{2} \io e^{v(\cdot, t)}  |\nabla a(\cdot, t)|^2\\
 &\le& \frac{1}{2} \io e^{v_0} |\nabla a_0|^2 + c_1\int_0^t \io v^2
       + c_2\int_0^t \io z^2 + \mu^2 c_3\int_0^t \io (1-ae^v)^2\\
 &\le& \frac{1}{2} \io e^{v_0} |\nabla a_0|^2 + c_1\int_0^\infty \io v^2
       + c_2\int_0^\infty \io z^2 + \mu^2 c_3\int_0^\infty \io (1-ae^v)^2\\
 &\le& \frac{1}{2} \io e^{v_0} |\nabla a_0|^2 + c_1c_4 + c_2 c_5
       + \mu^2 c_3 c_6
       \qquad\mbox{for all } t>0
 \eas
due to (\ref{13.4}) and (\ref{13.5}), which immediately entails
(\ref{13.1}) and (\ref{13.2}) because $e^v\ge 1$.
 \qed
Now using Lemma \ref{lem12} in conjunction with the compactness
properties implied by Lemma \ref{lem13} we can claim the following
$L^2$ stabilization feature of $a$.
\begin{lem}\label{lem14}
If $\beta<1$, then we have
  \be{14.1}
  \| a(\cdot, t)-1\|_{L^2(\Omega)} \to 0
  \qquad\mbox{as } t\to\infty.
  \ee
\end{lem}
\proof Combining Lemma \ref{lem12} with Lemma \ref{lem8} and
Corollary \ref{cor11} we readily see that
 \bas
 \int_0^\infty \io (a-1)^2<\infty.
 \eas
 Relying on this basic decay information together with the relative compactness properties implied by (\ref{13.1})
 and (\ref{13.2}), we use a straightforward argument by
 contradiction to achieve (\ref{14.1}); we may refer to \cite[Lemma 3.5]{taowin_262} for details in a quite similar
 setting.
 \qed
A simple interpolation yields a slightly stronger convergence
statement in the sense:
\begin{lem}\label{lem15}
 Let $\beta<1$. Then the solution of (\ref{0}) has the property
  \be{15.1}
    u(\cdot,t)  \to 1
    \quad \mbox{in } L^p(\Omega)
    \quad \mbox{for all } p\in (2,\infty)
    \qquad \mbox{as } t\to\infty.
  \ee
\end{lem}
\proof Since
 \bas
 \hspace*{-4mm}
 \|u-1\|_{L^p(\Omega)}
 &=& \|u-1\|_{L^\infty(\Omega)}^\frac{p-2}{p}
     \|u-1\|_{L^2(\Omega)}^\frac{2}{p}\\
 &=& \|e^va-1\|_{L^\infty(\Omega)}^\frac{p-2}{p}
     \|e^va-1\|_{L^2(\Omega)}^\frac{2}{p}\\
 &\le& \bigg\{e^{\|v_0\|_{L^\infty(\Omega)}}
        \|a\|_{L^\infty(\Omega)} +1\bigg\}^\frac{p-2}{p} \cdot
      \bigg\{\|a-1\|_{L^2(\Omega)} +
      \|e^va(1-e^{-v})\|_{L^2(\Omega)}\bigg\}^\frac{2}{p}\\
 &\le& \bigg\{e^{\|v_0\|_{L^\infty(\Omega)}}
        \|a\|_{L^\infty(\Omega)} +1\bigg\}^\frac{p-2}{p} \cdot
      \bigg\{\|a-1\|_{L^2(\Omega)} +
       e^{\|v_0\|_{L^\infty(\Omega)}}
        \|a\|_{L^\infty(\Omega)}\|v\|_{L^2(\Omega)}\bigg\}^\frac{2}{p}
 \eas
for all $t>0$ due to the H\"{o}lder inequality, (\ref{a}),
(\ref{vinfty}) and the basic inequality $1-e^{-s}\le s$ for all
$s\ge 0$, this straightly results from Lemma \ref{lem14}, Corollary
\ref{cor11} and Lemma \ref{lem8}. \qed
Up to now, we have achieved all goals in Theorem \ref{theo1}.\abs
\proofc of Theorem \ref{theo1}. \quad
  This immediately is the outcomes of Lemma \ref{lem_basic}, Lemma
  \ref{lem8}, Lemma \ref{lem9}, Corollary \ref{cor11} and Lemma \ref{lem15}.
\qed

\end{document}